\documentstyle[12pt]{article}
\newtheorem{lemma}{Lemma}
\newtheorem{theorem}{Theorem}
\newtheorem{prop}{Proposition}[theorem]
\newtheorem{corollary}{Corollary}

\let\Bbb=\bf

\newcommand{\newsection}{
\section}
\def\appendix#1{
\addtocounter{section}{1} \setcounter{equation}{0}
\renewcommand{\thesection}{\Alph{section}}
\section*{Appendix \thesection\protect\indent
#1}
}

\def\be{\begin{equation}}
\def\ee{\end{equation}}
\def\bea{\begin{eqnarray}}
\def\eea{\end{eqnarray}}

\let\Bbb=\bf
\def\HH{{\Bbb H}}
\def\RR{{\Bbb R}}
\def\CC{{\Bbb C}}
\def\ZZ{{\Bbb Z}}

\def\Li{\mathop{\mbox{Li}}\nolimits}

\begin{document}
\setlength{\unitlength}{1.5mm}

\begin{center}
{\Large Quantum Teichm\"uller space}\\
\vspace{4mm}
\end{center}
\begin{center}
\vspace{4mm}
L. Chekhov\\
{\it Steklov Mathematical Institute,}\\
{\it  Gubkina 8, 117966, GSP--1, Moscow, Russia} \\
{\it E-mail: {\tt chekhov@mi.ras.ru}}\\
and \\
V.V. Fock \\
{\it Institute of Theoretical and Experimental Physics,}\\
{\it B. Cheremushkinskaya 25, 117259 Moscow, Russia}\\
{\it E-mail: {\tt fock@heron.itep.ru}}
\vspace{18pt}
\end{center}
\vskip 0.9 cm
\begin{abstract}
We describe explicitly a noncommutative deformation of the $*$-algebra of
functions on the Teichm\"uller space of Riemann surfaces with holes
equivariant w.r.t. the mapping class group action.
\end{abstract}

\renewcommand{\thefootnote}{\arabic{footnote}}
\setcounter{footnote}{0}

\newsection{Introduction}
The main goal of the present paper is to describe the Hilbert space and the
space
of observables of quantum gravity in 2+1 dimensions. For this
purpose, we
use the following scheme. According to E.Verlinde and H.Verlinde~\cite{VV}
and Witten~\cite{Witten} (see also the monograph by Carlip~\cite{C}),
the classical phase space of Einstein gravity on a 3D manifold is the
Teichm\"uller space of its 2D boundary. (This
is analogous to the fact that
the classical phase space for 3D Chern--Simons theory is the moduli space of
flat connections on the boundary manifold.)
The Teichm\"uller space possesses the
canonical (Weil--Petersson) Poisson structure and the mapping class group as
a symmetry group. According to the correspondence principle ~(1) the algebra
of observables of the corresponding quantum theory is a noncommutative
deformation of the $*$-algebra of functions on the classical phase space in
the direction given by the Poisson structure; (2) the Hilbert spaces of the
theory is the $*$-representation space of this algebra; (3) the symmetry
group acts on the algebra of observables by automorphisms. To solve this
quantization
problem it suffices to construct a family of $*$-algebras, which depends
on the quantization parameter~$\hbar$, determine
the action of the mapping class
group by outer automorphisms, and show that the  algebra and the action
reproduces the classical algebra, the classical action, and the Poisson
structure in the limit $\hbar \rightarrow 0$ (provided we believe in the
existence and uniqueness of the quantization of a Poisson manifold).

In the present paper, we solve a very closely related problem,
however not
exactly the problem described above. We consider open 2D surfaces
(although this case can be also interpreted in the spirit of
the (2+1)D gravity, we skip the
discussion to avoid long definitions). The corresponding Teichm\"uller space
has a degenerate Weil--Petersson Poisson structure and the mapping class
group is a symmetry group. We describe the deformation quantization of the
corresponding Teichm\" uller space, the action of the mapping class group by
outer automorphisms, the unitary representations of the algebra, and
the induced action of the mapping class group on the representation space.

We show that the action of the mapping class group on the algebra of
observables depends on the quantization parameter~$\hbar$ but
remains unchanged if we replace~$\hbar$ by $1/\hbar$, which is in
agreement with a similar symmetry observed in \cite{ZZ}.

Note that according to~\cite{VV}, the representation space of the
observable algebra can be also interpreted as the space of conformal blocks
of the Liouville conformal field theory. The construction described in this
paper can be interpreted therefore as the construction of certain conformal
block spaces and the mapping class group actions for this
conformal field theory.

Briefly, the structure of the paper is the following. We give a very simple
description of the mapping class group in terms of generators and relations
of a certain groupoid called the {\em modular groupoid} (Sec.~2). We
introduce coordinates on the Teichm\"uller space, such that the mapping
class group action and the Weil--Petersson form takes an
especially simple form simple
(Sec.~3).  Then, we deform the algebra of functions on the Teichm\"uller
space into a Weyl algebra, define the action of the generators of the modular
groupoid, and verify that the relations between the generators are indeed
satisfied (Sec.~4).

  The main mathematical ingredient of the construction is a version of the
quantum
logarithm and dilogarithm by Faddeev and Kashaev~\cite{FK}. We
interpret the corresponding five-term relation as the only nontrivial
relation in the modular groupoid. The properties of the
quntum logarithm and dilogarithm functions are proven in Sec.~4.1.

  A similar construction was done independently and simultaneously by
Kashaev~\cite{kashaev}. However, our construction seems to be more
elementary and universal.


\newsection{Mapping class group. Graph description.}
Recall that the mapping class group ${\cal D}(S)$ of a 2D surface $S$ is the
group of homotopy classes of diffeomorphisms of the surface $S$. In
this section, we give a simple combinatorial description ${\cal D}(S)$ for
any open surface $S$.

A {\it fat graph} is a graph with a given cyclic ordering of
ends of edges entering each vertex. Any fat graph have the distinguished set
of closed paths called {\em faces}. These are the paths turning left at each
vertex. Gluing annuli to a graph along its faces produces a 2D surface with
the graph embedded into it. Therefore, we have a correspondence between the
fat graphs and topological types of open 2D surfaces. Any open surface
corresponds to a nonempty finite set of graphs (these are exactly the
graphs that can be embedded into the surface in such a way that the surface
can be shrunken onto it and the cyclic order of ends of edges coincides with
the one induced by the orientation of the surface).

Any fat graph has its dual --- the fat graph with
the vertices of the dual graph
being the faces of the original one. Each edge of the original graph
corresponds to the dual graph edge that
connects the faces the initial edge is incident to.

Denote by $|\Gamma|(S)$ the set of combinatorial types of three-valent
graphs corresponding to a given surface. For any element of $|\Gamma|(S)$
fix a {\em marking}, i.e., a numeration of the edges. Denote by
$\Gamma(S)$ the set of isotopy classes of embeddings of marked fat graphs
into $S$. The presence of the marking changes the set of embeddings since
some graphs may have nontrivial symmetry group. Introducing the marking is a
tool to remove this symmetry. (Here and below the vertical lines $|\cdot|$
indicate the diffeomorphism class.)

The mapping class group ${\cal D}(S)$ obviously acts freely on the space of
embedded marked graphs having the space of combinatorial graphs as a
quotient,
$$
\Gamma(S)/{\cal D}(S) = |\Gamma|(S).
$$

Recall that a group can be thought
of to be a category with only one object
and with all morphisms being invertible. Analogously, a {\em groupoid} is just
a category such that all morphisms are invertible and such that any two
objects are related by at least one morphism. Since the
automorphism groups of
different objects of a groupoid are obviously isomorphic to
each other, we can associate a group to a groupoid in the canonical way. We
are going to construct the groupoid giving the mapping class group and
admitting a simpler description in terms of generators and relations than
the mapping class group itself.

\noindent{\bf Definition.}
Let the set $|\Gamma|(S)$ be the set of objects. For any two graphs
$|\Gamma|, |\Gamma_1| \in |\Gamma(S)|$ let a morphism from $|\Gamma|$ to
$|\Gamma_1|$ be a homotopy class of marked embeddings of both $|\Gamma|$ and
$|\Gamma_1|$ into $S$ modulo the diagonal mapping class group action; we
denote this morphism
by $|\Gamma,\Gamma_1|$. If we have three embedded marked graphs
$\Gamma,\Gamma_1,\Gamma_2$, then by definition the composition of
$|\Gamma,\Gamma_1|$ and $|\Gamma_1,\Gamma_2|$ is $|\Gamma,\Gamma_2|$. The
above described category is called the {\em
modular groupoid}.

One can easily verify that ({\rm 1}) the multiplication of morphisms is
unambiguously defined; ({\rm 2}) the class of the diagonal embedding
$|\Gamma,\Gamma|$ is the identity morphism and the inverse of the morphism
$|\Gamma,\Gamma_1|$ is $|\Gamma_1,\Gamma|$; ({\rm 3}) the group of
automorphisms of an object is the mapping class group ${\cal D}(S)$.

To give a description of the modular groupoid by generators and relations we
need to introduce the distinguished sets of morphisms called {\em flips} and
{\em graph symmetries}. We call a morphism $|\Gamma,\Gamma_\alpha|$ a flip
if the embedding $\Gamma_\alpha$ is obtained from the embedding $\Gamma$ by
shrinking an edge $\alpha$ and blowing up the obtained four-valent vertex in
the other
direction (see Fig.~3 below). We use the notation $\Gamma_\alpha$ in order to
emphasize the relation of this graph to the graph $\Gamma$.
Note that for the given graph~$\Gamma$,
several marked embedded graphs may be denoted by $\Gamma_\alpha$
because no marking of
$\Gamma_\alpha$ is indicated.

To each symmetry $\sigma$ of a graph $\Gamma$ we associate an
automorphism, which is just $|\Gamma, \Gamma \sigma|$.

There is no canonical identification of edges of different
graphs even if a morphism between them is given.  However, for two graphs
related by a flip, we can introduce such
an identification. It is especially
transparent in the dual picture where a flip just replaces one edge by
another. Hence, one can identify the set of edges of two graphs as far as a
representation of a morphism between the graphs as a sequence of flips is
given. We exploit this identification and denote the corresponding edges
of different graphs by the same letter if it is clear
which sequence of flips relating
these graphs is considered. To avoid confusion, note that this
identification has nothing to do with the marking.

In this notation, the graph $\Gamma_{\alpha_1\cdots\alpha_n}$ is the graph
obtained as a result of consecutive flips
$\alpha_n$, \ldots, $\alpha_1$ of edges of a given graph~$\Gamma$.

There are three kinds of relations between flips, which are satisfied for
any choice of marking for the graphs entering the relations.

\begin{prop}

A square of a flip is a graph symmetry\/{\rm:}
if $|\Gamma_\alpha
,\Gamma|$ is a flip in an edge $\alpha$, then $|\Gamma, \Gamma_\alpha|$
is also a flip and\footnote{The notation {\bf R.n} indicates the number
$n$ of graphs entering this relation.}

\noindent{\bf R.2.}\ \ \
$|\Gamma, \Gamma_\alpha| |\Gamma_\alpha, \Gamma|=1$.

Flips in disjoint edges commute\/{\rm:} if
$\alpha$ and $\beta$
are two edges having no common vertices, then

\noindent{\bf R.4.} \ \ \
$|\Gamma_{\alpha\beta}, \Gamma_\alpha||\Gamma_\alpha, \Gamma| =
|\Gamma_{\alpha\beta},
\Gamma_\beta||\Gamma_\beta, \Gamma|$.

Five consecutive flips in edges $\alpha$ and $\beta$
having one common vertex is the identity\/{\rm:}
for such~$\alpha$
and~$\beta$, the graphs $\Gamma_{\alpha\beta}$ and $\Gamma_{\beta\alpha}$ are
related by a flip and

\noindent{\bf R.5.} $|\Gamma,\Gamma_\alpha||\Gamma_\alpha,
\Gamma_{\beta\alpha}||\Gamma_{\beta\alpha},\Gamma_{\alpha\beta}||
\Gamma_{\alpha\beta},\Gamma_\beta||\Gamma_\beta,\Gamma| =1$.
\end{prop}

The proofs of relations {\bf R.2} and {\bf R.4}
are obvious. Relation {\bf R.5} can be seen more
transparently in the dual graph picture.
Indeed, a graph dual to a three-valent
graph is a graph having triangular faces. A flip of the original graph
corresponds to removing an edge on the dual graph and inserting another
diagonal of the appearing quadrilateral. Figure~1 shows that
the combination of the five flips is the identity. \hfill{\bf q.e.d.}

\centerline{\unitlength 0.45mm
\begin{picture}(0,0)(100,-145)
\put(0,0){\line(2,1){ 40}}
\put(40,20){\line(2,-1){ 40}}
\put(80,0){\line(-1,-2){ 20}}
\put(60,-40){\line(-1,0){40}}
\put(20,-40){\line(-1,2){ 20}}
\put(0,0){\line(1,0){80}}
\put(0,0){\line(3,-2){60}}
\put(20,13){\makebox(0,0)[rb]{$E_0$}}
\put(60,13){\makebox(0,0)[lb]{$D_0$}}
\put(8,-21){\makebox(0,0)[rc]{$A_0$}}
\put(72,-21){\makebox(0,0)[lc]{$C_0$}}
\put(40,-43){\makebox(0,0)[ct]{$B_0$}}
\put(40,3){\makebox(0,0)[cb]{$Y_0$}}
\put(29,-22){\makebox(0,0)[ct]{$X_0$}}
\put(45,-50){\vector(1,-2){5}}
\end{picture}
\begin{picture}(0,0)(20,-200)
\put(0,0){\line(2,1){ 40}}
\put(40,20){\line(2,-1){ 40}}
\put(80,0){\line(-1,-2){ 20}}
\put(60,-40){\line(-1,0){40}}
\put(20,-40){\line(-1,2){ 20}}
\put(0,0){\line(3,-2){60}}
\put(40,20){\line(1,-3){20}}
\put(20,13){\makebox(0,0)[rb]{$C_4$}}
\put(60,13){\makebox(0,0)[lb]{$B_4$}}
\put(8,-21){\makebox(0,0)[rc]{$D_4$}}
\put(72,-21){\makebox(0,0)[lc]{$A_4$}}
\put(40,-43){\makebox(0,0)[ct]{$E_4$}}
\put(53,-8){\makebox(0,0)[lc]{$X_4$}}
\put(29,-22){\makebox(0,0)[ct]{$Y_4$}}
\put(1,-29){\vector(-1,-1){10}}
\end{picture}
\begin{picture}(0,0)(-60,-145)
\put(0,0){\line(2,1){ 40}}
\put(40,20){\line(2,-1){ 40}}
\put(80,0){\line(-1,-2){ 20}}
\put(60,-40){\line(-1,0){40}}
\put(20,-40){\line(-1,2){ 20}}
\put(40,20){\line(1,-3){20}}
\put(40,20){\line(-1,-3){20}}
\put(20,13){\makebox(0,0)[rb]{$A_3$}}
\put(60,13){\makebox(0,0)[lb]{$E_3$}}
\put(8,-21){\makebox(0,0)[rc]{$B_3$}}
\put(72,-21){\makebox(0,0)[lc]{$D_3$}}
\put(40,-43){\makebox(0,0)[ct]{$C_3$}}
\put(53,-8){\makebox(0,0)[lc]{$Y_3$}}
\put(27,-8){\makebox(0,0)[rc]{$X_3$}}
\put(13,21){\vector(-1,1){10}}
\end{picture}
\begin{picture}(0,0)(-30,-65)
\put(0,0){\line(2,1){ 40}}
\put(40,20){\line(2,-1){ 40}}
\put(80,0){\line(-1,-2){ 20}}
\put(60,-40){\line(-1,0){40}}
\put(20,-40){\line(-1,2){ 20}}
\put(40,20){\line(-1,-3){20}}
\put(20,-40){\line(3,2){60}}
\put(20,13){\makebox(0,0)[rb]{$D_2$}}
\put(60,13){\makebox(0,0)[lb]{$C_2$}}
\put(8,-21){\makebox(0,0)[rc]{$E_2$}}
\put(72,-21){\makebox(0,0)[lc]{$B_2$}}
\put(40,-43){\makebox(0,0)[ct]{$A_2$}}
\put(27,-8){\makebox(0,0)[rc]{$Y_2$}}
\put(52,-23){\makebox(0,0)[ct]{$X_2$}}
\put(58,19){\vector(1,2){6}}
\end{picture}
\begin{picture}(0,230)(70,-65)
\put(0,0){\line(2,1){ 40}}
\put(40,20){\line(2,-1){ 40}}
\put(80,0){\line(-1,-2){ 20}}
\put(60,-40){\line(-1,0){40}}
\put(20,-40){\line(-1,2){ 20}}
\put(20,-40){\line(3,2){60}}
\put(0,0){\line(1,0){ 80}}
\put(20,13){\makebox(0,0)[rb]{$B_1$}}
\put(60,13){\makebox(0,0)[lb]{$A_1$}}
\put(8,-21){\makebox(0,0)[rc]{$C_1$}}
\put(72,-21){\makebox(0,0)[lc]{$E_1$}}
\put(40,-43){\makebox(0,0)[ct]{$D_1$}}
\put(52,-23){\makebox(0,0)[ct]{$Y_1$}}
\put(40,3){\makebox(0,0)[cb]{$X_1$}}
\put(82,-21){\vector(1,0){15}}
\end{picture}}
\vspace{5mm}
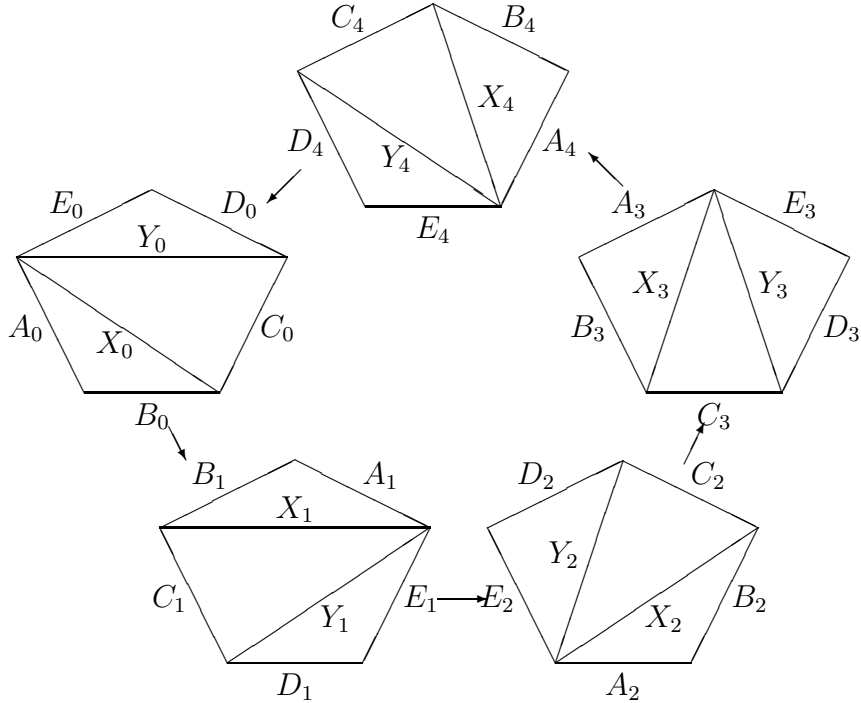
\centerline{\bf Fig.~1}
\vspace{5mm}

\begin{theorem}~ \label{th1}

\noindent{\bf 1.} Flips and graph symmetries generate the modular groupoid.

\noindent{\bf 2.} The only relations between the generators are {\bf R.2},
{\bf R.4}, {\bf R.5}, and the natural relations between flips and graph
symmetries.
\end{theorem}

Replacing the mapping class group by the modular groupoid,
we can simply express the latter through generators and relations.

Note that a graph symmetry can be represented as a ratio of two flips in
a given edge and the modular groupoid is therefore generated
only by the flips.
We do not describe relations between flips and graph symmetries in
details because they are quite obvious. In fact, the symmetry groups of
$\Gamma$ and $\Gamma_\alpha$ act transitively on on the set of flips
$|\Gamma, \Gamma_\alpha|$, and this action can be considered as relations
between flips and graph symmetries.

Theorem~1 can be proved using direct combinatorial methods of the simplicial
geometry (cf. Viro \cite{Viro}). However, we give the
main idea of another proof, which is
more specific for the 2D situation.

\noindent{\bf Proof of Theorem~1.}
To any connected simplicial complex ${\cal S}$ we can
associate a groupoid by taking a point in each top-dimensional simplex for
objects and  the homotopy classes of oriented paths connecting the chosen
points as morphisms. The corresponding group is the fundamental group of the
topological space given by the complex.

To any codimension one simplex we can associate two classes of paths
(differing by orientation and having the identity morphism
as their product)
connecting adjacent top-dimensional simplices. It is natural to
call them {\em flips}. We can associate a
relation between the flips to any codimension two simplex.
It is obvious that this set of flips generates
the groupoid and that the only relation between the flips are given by
codimension two simplices.

The same is true for an orbifold simplicial complex, where we replace
simplices by quotients of simplices by finite groups. In this case, we must
choose one generic point per each top-dimensional simplex as an object and
orbifold homotopy classes of paths as morphisms. The corresponding group is
the orbifold fundamental group of the orbifold given by the complex. The
groupoid is now generated by flips and groups of top dimension simplices and
still the only nontrivial relations are those given by codimension two
simplices.

Consider now the Strebel \cite{Streb1} orbifold simplicial decomposition of
the moduli space of complex structures on $S$. The orbifold fundamental
group of the moduli space ${\cal M}$ is just the mapping class group ${\cal
D}(S)$. Recall that
Strebel orbisimplices are enumerated by fat graphs corresponding to~$S$
and the
dimension of a simplex is equal to the number of its edges.  One can
easily see that the groupoid of the Strebel complex coincides with the
modular
groupoid. Moreover, the flips of the former correspond to the flips of the
latter.  The relations between flips are given by codimension two cells,
which correspond either to graphs with two four-valent vertices
(which produces
relation {\bf R.4}) or to graphs with one five-valent vertex (which produces
relation {\bf R.5}). Relation {\bf R.2} holds true for any
symplicial complex.\hfill{\bf q.e.d.}

\section{Teichm\"uller space. Graph description.}

Recall briefly the combinatorial description of the Teichm\"uller spaces of
complex structures on Riemann surfaces with holes.
The Teichm\"uller space ${\cal T}^h(S)$ is the space of complex structures on
a (possibly open) surface $S$ modulo diffeomorphisms homotopy equivalent to
the identity.
There are two types of behaviour of a complex structure in a
vicinity of a hole. The vicinity can be isomorphic (as a complex manifold)
either to an annulus or to a punctured disc. The holes of the second kind
are called punctures.

For technical reasons, instead of the Teichm\"uller space ${\cal T}^h(S)$, we
consider its finite covering ${\cal T}^H(S)$. A point of ${\cal
T}^H(S)$ is determined by a point of ${\cal T}^h(S)$ and by orientations of 
all
holes of $S$ which are not punctures. This covering is obviously ramified
over the subspace of surfaces with punctures.

The Poincar\'e uniformization theorem states that any complex surface $S$ is
a quotient of the upper half-plane $\HH$ under the action of a discrete
(Fuchsian) subgroup $\Delta(S)$ of the group $PSL(2,\RR)$ of
the automorphisms of~$\HH$,
$$
S={\HH}/\Delta(S).
$$
The upper half plane possesses the $PSL(2,\RR)$-invariant Hermitian metric
with curvature -1 given in the standard coordinates $z, \overline{z}$ by
$(\mbox{Im~}z)^{-2}dz d\overline{z}$. Hence on any complex surface $S$ there
exists a canonical Hermitian metric of the curvature $-1$.

Any homotopy class of closed curves $\gamma$ in $S$
(except the curves surrounding one
puncture) contains a unique closed geodesic of the length
\be
\label{length}
l(\gamma)=\left|\log
{\lambda_1}/{\lambda_2}\right|,
\ee
where $\lambda_1$ and $\lambda_2$ are different
eigenvalues of the element of $PSL(2,\RR)$ that corresponds to
${\gamma}$. (If we assume that geodesics' surrounding punctures have zero
length, then formula (\ref{length}) is valid for all classes
of curves.)

Since the work of Penner \cite{Penner}, the fat graphs are used for
describing not only moduli, but also Teichm\"uller spaces. We use a
version of this description, which is rather explicit and simple.

\begin{theorem}
Given a three-valent marked embedded graph $\Gamma \in
\Gamma(S)$, there exists a canonical isomorphism between the set of points of
the Teichm\"uller space ${\cal T}^H(S)$ and the set ${\bf R}^{\hbox{\rm\small
\#\,edges}}$ of assignments of real numbers to the edges of the graph.
\end{theorem}

For details of the proof see \cite{Fock2}. The construction of the numbers
on edges starting from a point of ${\cal T}^H(S)$ and {\em vice versa} is
rather explicit and elementary. As an illustration, we
construct the
Fuchsian group $\Delta(S)\subset PSL(2,\RR)$ corresponding to a given set of
numbers on edges of a graph $\Gamma \in |\Gamma|(S)$. To describe the
Fuchsian group, we must associate an element $X_\gamma\in PSL(2,\RR)$ to
any element of the fundamental group $\gamma \in \pi_1(S)$.

Starting from a fat graph $\Gamma \in \Gamma(S)$, consider another fat graph
$\widetilde{\Gamma}$ in which
a vicinity of each vertex is removed and the arising
three ends of edges are connected by three
new edges forming a triangle. Orient
all edges of the inserted triangles clockwise as in Fig.~2.

\vspace{10pt}
\setlength{\unitlength}{1mm}%
\begin{picture}(50,27)(-50,48)
\thinlines
\put(28,62){\vector( 1, 0){ 5}}
\thicklines
\put(38,70){\line( 1,-2){ 2}}
\put(48,62){\line( 1, 0){16}}
\put(72,66){\line( 1, 2){ 2}}
\put(72,58){\line( 1,-2){ 2}}
\put(38,54){\line( 1, 2){ 2}}
\put(40,58){\vector( 0, 1){ 8}}
\put(48,62){\vector( -2, -1){ 8}}
\put(40,66){\vector( 2, -1){ 8}}
\put(72,66){\vector( 0, -1){ 8}}
\put(72,58){\vector( -2, 1){ 8}}
\put(64,62){\vector( 2, 1){ 8}}
\thicklines
\put(-12,70){\line( 1,-2){ 4}}
\put(-8,62){\line( 1, 0){28}}
\put(20,62){\line( 1, 2){ 4}}
\put(20,62){\line( 1,-2){ 4}}
\put(-8,62){\line(-1,-2){ 4}}
\end{picture}

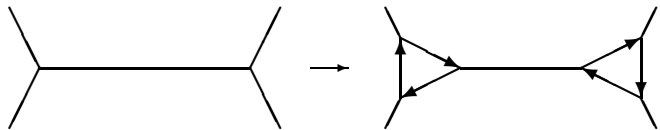
\centerline{\bf Fig.~2}
\vspace{5mm}

We now introduce the matrices $X(z), z \in \RR$ and $L$ from
$PSL(2,\RR)$,
$$
\begin{array}{rcl}
X(z)&=&\left(\begin{array}{cc}
                 0&e^{z/2}\\
                 -e^{-z/2} & 0
                     \end{array}\right),\\
L&=&\left(\begin{array}{cc}
                1& 1\\
                 -1&0
                     \end{array}\right),
\end{array}
$$

To any path $\gamma$ on the graph $\widetilde{\Gamma}$
consisting of $N$ edges we
associate an element of $X_\gamma \in PSL(2,\RR)$ by the rule $X_\gamma =
X_N\ldots X_1$, where $X_i = X(z_{\alpha})$, 
$L$ or $L^{-1}$ depending on whether the
$i$th segment of the path goes along the long edge 
or along a side of a small triangle in
positive or negative direction, respectively.

Choosing one vertex and considering all elements of $PSL(2,\RR)$
corresponding to paths starting and ending at this vertex, we
obtain the desired Fuchsian group $\Delta(S)$.

Because $X(z)^2 = 1$ in $PSL(2,\RR)$ we indeed obtain a homomorphism of
$\pi_1(\tilde{\Gamma}) \rightarrow PSL(2,\RR)$. Note also that since
$L^3=1$, this homomorphism factorises as
$\pi_1(\widetilde{\Gamma}) \rightarrow \pi_1(\Gamma)
\rightarrow PSL(2,\RR)$ and
the constructed group $\Delta(S)$ is in fact isomorphic to $\pi_1(\Gamma) =
\pi_1(S)$ (in other words, the product along $\gamma$ does not depend on how
$\gamma$ goes aroung small triangles).

Using this construction of the Fuchsian group, we can compute the length of
any closed geodesic using (\ref{length}). In particular, the length
of the geodesics surrounding a hole (face)~$f$ is
\be
l_f = \left|\sum_{\alpha \in f} z_\alpha\right|. \label{lface}
\ee
The sign of the sum determines the mutual orientation of a hole and
the surface. (For a generic path, the length
is obviously nonlinear in $z_\alpha$'s.)

\vspace{10pt}
\setlength{\unitlength}{1.5mm}%
\begin{picture}(50,27)(-12,48)
\thicklines
\put(28,70){\line( 1,-2){ 4}}
\put(32,62){\line( 1, 0){28}}
\put(60,62){\line( 1, 2){ 4}}
\put(60,62){\line( 1,-2){ 4}}
\put(32,62){\line(-1,-2){ 4}}
\thinlines
\put(18,62){\vector(-1, 0){  0}}
\put(18,62){\vector( 1, 0){ 5}}
\thicklines
\put(10,54){\line( 2,-1){ 8}}
\put(10,70){\line( 0,-1){16}}
\put(10,54){\line(-2,-1){ 8}}
\put( 2,74){\line( 2,-1){ 8}}
\put(10,70){\line( 2, 1){ 8}}
\put( 4,74){\makebox(0,0)[lb]{$a$}}
\put(16,74){\makebox(0,0)[rb]{$b$}}
\put(12,62){\makebox(0,0)[lc]{$z$}}
\put(16,50){\makebox(0,0)[rt]{$c$}}
\put( 4,50){\makebox(0,0)[lt]{$d$}}
\put(30,54){\makebox(0,0)[lt]{$d - \phi(-z)$}}
\put(62,54){\makebox(0,0)[rt]{$c+\phi(z)$}}
\put(62,69){\makebox(0,0)[rb]{$b-\phi(-z)$}}
\put(30,69){\makebox(0,0)[lb]{$a+\phi(z)$}}
\put(47,64){\makebox(0,0)[cb]{$-z$}}
\end{picture}

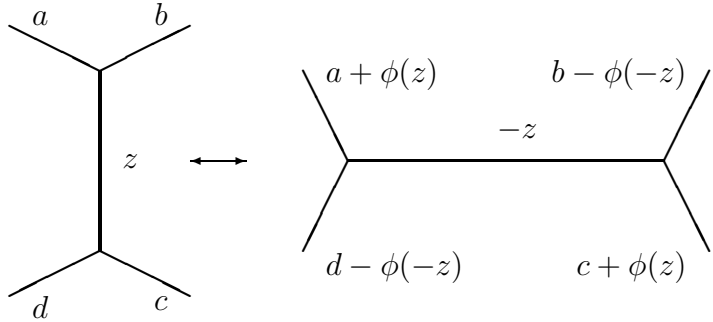
\centerline{\bf Fig.~3}
\vspace{5mm}

Therefore, for each embedded graph $\Gamma \in \Gamma(S)$, we obtain a
global coordinate system on ${\cal T}^H(S)$ and transition maps
between coordinate systems corresponding to different embedded graphs
must exist. In
other words, there exists a transition map for any morphism of the modular
groupoid. As follows from Theorem~\ref{th1}, any transition map can be
expressed as a composition of the transition maps corresponding to flips.
For a flip the transition map is given by the following rule: the numbers on
all edges except the five involved in the flip remain unchanged and the
numbers on the remaining
five edges are changed as shown in Fig.~3, while
\be
\phi(z) = \log(e^z + 1).
\ee

Although it follows directly
from the geometrical setting, it can be straightforwardly verified
that this transformation indeed does not change the numbers on edges when
the sequences of flips {\bf R.2}, {\bf R.4}, and {\bf R.5} are applied.

A canonical Poisson structure called the Weil--Peterson structure exists  on
${\cal T}^H(S)$ (cf. Goldman \cite{Goldman}). This structure is degenerate,
and the Casimir functions are just the lengths of geodesics surrounding
holes (\ref{lface}).

\begin{theorem} \label{th3}
In the graph coordinates, the Weil--Petersson bracket $B_{WP}$ has
a very simple form\/{\rm:}
\be
\label{WP-PB}
B_{WP} = \sum_{{\alpha}} \frac{\partial}{\partial
z_{{\alpha}}}\wedge \frac{\partial}{\partial z_{{\alpha}^r}},
\ee
where the sum is taken over all {\em oriented} edges ${\alpha}$, while
${\alpha}^r$ is the edge next to the right with respect to the end of
${\alpha}$ and $z_{{\alpha}}$ is the number assigned to the
{\rm(}nonoriented\/{\rm)} edge ${\alpha}$.
\end{theorem}

The idea of the proof was given in \cite{Fock2}. In fact, for
our purposes it is sufficient to check that
the form (\ref{WP-PB}) and the
set of Casimir functions are independent on the choice of the
embedded graph, what can be easily
verified explicitly.

\section{Quantization}
Recall that a quantization of a Poisson manifold equivariant w.r.t. a
discrete group action is a family of $*$-algebras ${\cal A}^\hbar$ depending
smoothly on a positive real parameter $\hbar$, acting on~$G$ by outer
automorphisms and
having the following relation to the Poisson manifold.

{\bf 1.} For $\hbar=0$, the algebra is isomorphic as a $G$-module to the
$*$-algebra of complex-valued function on the Poisson manifold.

{\bf 2.} The Poisson bracket on ${\cal A}^0$ given by $\{a_1, a_2\} =
\lim_{\hbar \rightarrow 0}\frac{[a_1,a_2]}{\hbar}$ coincides with the one
generated by the Poisson structure of the manifold.

We now quantize a Teichm\"uller space
${\cal T}^\hbar(S)$ equivariantly w.r.t.\ the mapping class group ${\cal
D}(S)$.

{\bf Remark.} By a smooth one-parameter family of $*$-algebras we mean just
that all these algebras can be identified as vector spaces and the
multiplication rule varies smoothly under changing
the parameter $\hbar$. We do not
require the identification
between algebras to be equivariant with respect to the group action.

Let ${\cal T}^\hbar(\Gamma)$ where $\Gamma \in |\Gamma|(S)$ is a
$*$-algebra generated by the generator $Z_\alpha^\hbar$
(one generator per one unoriented edge $\alpha$)
and relations
\be
[Z^\hbar_\alpha, Z^\hbar_\beta ] = 2\pi i\hbar\{z_\alpha, z_\beta\}
\label{qq}
\ee
with the $*$-structure
\be
(Z^\hbar_\alpha)^*=Z^\hbar_\alpha.
\ee
Here $z_\alpha$  and $\{\cdot,\cdot\}$ stand for the respective
coordinate functions on
the classical Teichm\"uller space and the Weil--Petersson Poisson bracket on
it. Note that according to formula (\ref{WP-PB}), the
right-hand side of (\ref{qq})
is merely a constant which may take only five values $0$, $\pm
2\pi i \hbar$, $\pm 4 \pi i \hbar$.

We have constructed one $*$-algebra per object of the modular groupoid. Now
in order to describe the equvariance, we must associate a
homomorphism of the corresponding $*$-algebras to any morphism of the
modular groupoid.
For this, we must associate a  unitary morphism
of $*$-algebras to any flip and any graph symmetry and verify that the
relations {\bf R.2}, {\bf R.4}, and {\bf R.5} are satisfied by these
morphisms.

The action of a graph symmetry can be obviously reduced to
permutating the generators.

We now define the flip morphisms by the rule
\be
\{A,B,C,D,Z\}\to\{A+\phi^\hbar(Z),B-\phi^\hbar(-Z),C+\phi^\hbar(Z),
D-\phi^\hbar(-Z),-Z\},
\label{q-mor}
\ee
where $A$, $B$, $C$, $D$, and~$Z$ are as in Fig.~4 and
$\phi^\hbar(x)$ is the real function of one real variable,
\begin{equation} \label{phi}
\phi^\hbar(z) =
-\frac{\pi\hbar}{2}\int_{\Omega} \frac{e^{-ipz}}{\sinh(\pi p)\sinh(\pi \hbar
p)}dp,
\end{equation}
while the contour $\Omega$ goes along the real axis bypassing the
singularity at the origin from above. (The generators on
the edges not shown on the Fig.~4 do not change)

\vspace{10pt}
\setlength{\unitlength}{1.5mm}%
\begin{picture}(50,27)(-12,48)
\thicklines
\put(28,70){\line( 1,-2){ 4}}
\put(32,62){\line( 1, 0){28}}
\put(60,62){\line( 1, 2){ 4}}
\put(60,62){\line( 1,-2){ 4}}
\put(32,62){\line(-1,-2){ 4}}
\thinlines
\put(18,62){\vector(-1, 0){  0}}
\put(18,62){\vector( 1, 0){ 5}}
\thicklines
\put(10,54){\line( 2,-1){ 8}}
\put(10,70){\line( 0,-1){16}}
\put(10,54){\line(-2,-1){ 8}}
\put( 2,74){\line( 2,-1){ 8}}
\put(10,70){\line( 2, 1){ 8}}
\put( 4,74){\makebox(0,0)[lb]{$A$}}
\put(16,74){\makebox(0,0)[rb]{$B$}}
\put(12,62){\makebox(0,0)[lc]{$Z$}}
\put(16,50){\makebox(0,0)[rt]{$C$}}
\put( 4,50){\makebox(0,0)[lt]{$D$}}
\put(30,54){\makebox(0,0)[lt]{$D - \phi^\hbar(-Z)$}}
\put(62,54){\makebox(0,0)[rt]{$C+\phi^\hbar(Z)$}}
\put(62,69){\makebox(0,0)[rb]{$B-\phi^\hbar(-Z)$}}
\put(30,69){\makebox(0,0)[lb]{$A+\phi^\hbar(Z)$}}
\put(47,64){\makebox(0,0)[cb]{$-Z$}}
\end{picture}

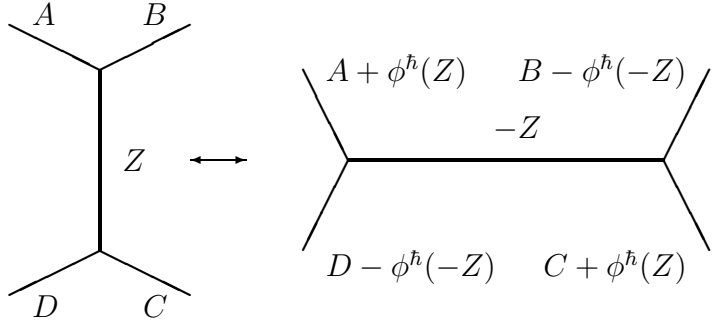
\centerline{\bf Fig.~4}
\vspace{5mm}

We now formulate the main theorem of the paper.

\begin{theorem} \label{th4}
The family of algebras ${\cal T}^\hbar(\Gamma)$ is a quantization of
the space ${\bf \cal T}^H(S)$. Relation (\ref{q-mor})
defines the action of the modular groupoid (and of
the mapping class group ${\cal D}(S)$, in particular) on ${\cal
T}^\hbar(\Gamma)$ by external $*$-morphisms.
\end{theorem}

Theorem \ref{th4} implies, in particular,
that the algebras corresponding to different
graphs are isomorphic, which is emphasized by
using the notation ${\cal T}^\hbar(S)$ instead of ${\cal T}^\hbar(\Gamma)$.

Before proving Theorem \ref{th4}, we discuss its
corollaries and the
properties of the constructed algebra and the mapping
class group action.

\begin{corollary}
{\bf C.1}. In the limit $\hbar \rightarrow 0$, the morphism
{\rm(\ref{q-mor})} coincides with the classical morphism in Fig~3.

\noindent{\bf C.2} The morphisms ${\cal T}^\hbar(\Gamma)\rightarrow {\cal
T}^{1/\hbar}(\Gamma)$ given by the mapping $Z^\hbar_\alpha
\mapsto Z^{1/\hbar}_\alpha/\hbar$ commute with the modular morphisms.

\noindent{\bf C.3}. The center of the algebra ${\cal T}^\hbar(S)$
is generated by the sums
$Z^\hbar_f =\sum_{\alpha \in f}{Z^\hbar_\alpha}$ ranging all edges $\alpha$
surrounding a given face $f$. The mapping class group acts on the generators
$Z^\hbar_f$ by permutations, i.e., exactly as
it acts on the holes of the surface $S$.

\noindent{\bf C.4}. For any set of real numbers $l_1, \ldots l_s$ assigned
to the holes of the surface $S$, there exists a unique unitary representation
$\pi^\hbar(S,l_1, \ldots l_s)$, of ${\cal T}^H(S)$ in a Hilbert space
$H^\hbar(S,l_1, \ldots l_s)$, such that the central elements are
represented by multiplications by these constants, $Z^\hbar_f = l_f$.

\noindent{\bf C.5}. Let ${\cal D}(S, l_1, \ldots, l_s) \subset {\cal D}(S)$
be a subgroup of the mapping class group preserving the assignment of the
numbers  $l_1, \ldots l_s$ to the holes. Then there exists a unique projective
representation $\rho$ of ${\cal D}(S, l_1, \ldots, l_s)$ in $H^\hbar(S,l_1,
\ldots, l_s)$, determined by the condition that $\rho(x) \pi(t) = \pi(xt)
\rho(x)$ for any $x \in {\cal D}(S, l_1, \ldots, l_s)$ and $t \in {\cal
T}^H(S)$.

\noindent{\bf C.6}. The obstacle for the representation $\pi(S, l_1, \ldots,
l_s)$ to be a true representation {\rm(}and not only a projective
one\/{\rm)} lies in
the second cohomology group $H^2({\cal D}(S, l_1, \ldots, l_s), U(1))$, or,
equivalently, in the second homology group of the orbifold moduli space
\newline $H^2({\cal T}^H(S)/{\cal D}(S, l_1, \ldots, l_s), U(1))$.
\end{corollary}

The main technical tool used in the proof of
Theorem \ref{th4} is a simple observation,
which we formulate as  a separate lemma.

\begin{lemma} \label{lem1}
Consider an algebra {\rm(}{\em quantum torus}{\rm)}
generated by two generators~$U$
and~$V$ satisfying the commutation relation $q UV = q^{-1} VU$, where
$q$ is a nonvanishing complex number. Then the automorphism of the algebra
given by $U \rightarrow (1+qU)V$, $V \rightarrow U^{-1}$ has order five.
\end{lemma}

We formulate few more corollaries of Theorem \ref{th4} and Lemma \ref{lem1}.
Although being technical, they seem
to be sufficiently elegant to be stated separately.

\begin{corollary}
{\bf C.7} \ Let $K$ be an operator acting in the Hilbert space
$L^2(\RR)$ and having the integral kernel
\be
K(x,z) = F^\hbar(z)e^{-\frac{zx}{2\pi i \hbar}},
\ee
where
\be
\label{dlc}
F^\hbar(z) = \exp\left(-\frac{1}{4}\int_\Omega \frac{e^{-ipz}}{p\sinh(\pi
p)\sinh(\pi \hbar p)}dp\right)
\ee
Then the operator $K$ is unitary up to a multiplicative constant and
satisfies the identity
\be
\label{pic}
K^5=\mbox{\em const}.
\ee

\noindent{\bf C.8} \ Let $\hbar=m/n$ be a rational number and assume that
both $m$ and $n$ are odd. Introduce a linear operator $L(u)$
acting in the space
$\CC^n$ and depending on one positive real parameter $u$
through its matrix
\be
\label{dld} L(u)^i_j = F^\hbar(j,u)q^{2ij},
\ee
where
$$
F^\hbar(j,u) =
(1+u)^{j/n}\prod_{k=0}^{j-1}(1+q^{2k-1}u^{1/n})^{-1}.
$$
Then the following
identity holds\/{\rm:}
\be
\label{pid}
L(v^{-1})L(u^{-1}v^{-1}+ u^{-1})L(v + v u^{-1}
+ u^{-1}) L(v + uv)L(u) = \mbox{\em const}.
\ee
\end{corollary}
\vspace{3mm}

\noindent {\bf Remarks.}

{\bf 1.} In accordance with the intuition, the set of irreducible
representations of ${\cal T}^\hbar(S)$ is in
the one-to-one correspondence with
the symplectic leaves of ${\cal T}^H(S)$.

{\bf 2.} The group ${\cal D}(S, l_1, \ldots l_s)$
actually depends only on the distribution
of values of the central generators $l_1, \ldots, l_s$. If
all $l_i$'s are equal, then, obviously,
${\cal D}(S, l_1, \ldots l_s) = {\cal D}(S)$.

{\bf 3.} The described constructions can be considered as
the generalization of
the Shale--Weil representation of the symplectic group (cf. \cite{Vergne}).
In the Shale--Weil case,
the Poisson manifold under quantization is the
first cohomology group $H^1(S)$ with compact
support and with the natural mapping class group action. This analogy becomes
even more transparent if one represents the space $H^1(S)$ as the space of
representations of $\pi_1(S)$ in $\RR$, and ${\cal T}^h(S)$ as an open
subset of the space of representations of $\pi_1(S)$ in $PSL(2,\RR)$.

{\bf 4.} There is a natural question whether the described set of mapping
class group representations constitutes a modular functor. In other words,
since for any embedding $S_1 \rightarrow S_2$ of one surface into another
there exists the induced homomorphism of the mapping class groups ${\cal
D}(S_1) \rightarrow {\cal D}(S_2)$, the question arises whether a constructed
unitary representation of ${\cal D}(S_2)$ always decomposes into a direct
sum of the constructed representations of ${\cal D}(S_1)$. The answer is
presumably negative.

{\bf 5.} The function $F^\hbar(z)$ (see (\ref{dlc})) asymptotically
behaves as $F(z) \rightarrow \exp\left\{\Li_2(e^z)/\hbar\right\}$ when $\hbar
\rightarrow 0$. Here $\Li_2(z)$ is the dilogarithm function
$$
\Li_2(e^z)
= \int_0^{e^z}\frac{\log(1+x)}{x}dx.
$$
Both this function and matrix
function (\ref{dld}) are versions of quantum dilogarithms introduced by
Faddeev and Kashaev~\cite{FK}. Note that the celebrated five-term relation
for $\Li_2$ can be rewritten in the form
$$
\widetilde{L}(v^{-1})
\widetilde{L}(u^{-1}v^{-1}+ u^{-1})
\widetilde{L}(v + v u^{-1} + u^{-1})
\widetilde{L}(v + uv)\widetilde{L}(u)=1,
$$
where $\widetilde{L}(u) = e^{\Li_2(u)}$.

{\bf 6.} Taking into account the previous remarks, the operators
representing the mapping class group in the Hilbert space $H^\hbar(S,l_1,
\ldots, l_s)$ can be written as compositions of ({\rm 1}) Fourier
transforms, ({\rm 2}) operators of multiplication by exponents of quadratic
forms, and ({\rm 3}) operators of the
multiplication by the function $F^\hbar$ (which is also
very close to the Shale--Weil representation, where only the first two
operations are involved).

{\bf 7.} Assuming the standard relation between
deformation quantization and geometric quantization to be true,
the algebra ${\cal T}^\hbar(S)$
can be represented in the space of sections of a certain linear
bundle over the Teichm\"uller space. Analogously, the subalgebra ${\cal
M}^\hbar(S) \subset {\cal T}^\hbar(S)$ that is
invariant under the action of the
mapping class group ${\cal D}(S)$ must act on sections of a certain linear
bundle over the moduli space. Modular forms are good candidates for these
sections.

\subsection{Properties of the function $\phi^h(z)$}

for proving Theorem \ref{th4} and its corollaries we present here the
properties of the function
$$
\phi^\hbar(z) =
-\frac{\pi\hbar}{2}\int_{\Omega} \frac{e^{-ipz}}{\sinh(\pi p)\sinh(\pi \hbar
p)}dp.
$$

\begin{prop}
For the function $\phi^{\hbar}(z)$, we have

\noindent{\bf P.1} \ \ \
$\lim_{\hbar \rightarrow 0}\phi^\hbar(z) = \log(e^z + 1)$;

\vspace{3mm}

\noindent{\bf P.2} \ \ \
$\phi^\hbar(z)-\phi^\hbar(-z)=z$;

\vspace{3mm}

\noindent{\bf P.3} \ \ \
$\overline{\phi^\hbar(z)} = \phi^\hbar(\overline{z})$;

\vspace{3mm}

\noindent{\bf P.4} \ \ \
$\frac{1}{\hbar}\phi^\hbar(z) = \phi^{1/\hbar}(z/\hbar)$;

\vspace{3mm}

\noindent{\bf P.5} \ \ \
$\phi^\hbar(z+i\pi\hbar)-\phi^\hbar(z-i\pi\hbar) = \frac{2\pi i
\hbar}{e^{-z}+1}$;

\vspace{3mm}

\noindent{\bf P.6} \ \ \
$\phi^\hbar(z+i\pi)-\phi^\hbar(z-i\pi) = \frac{2\pi
i}{e^{-z/\hbar}+1}$;

\vspace{3mm}

\noindent{\bf P.7} \ \ \
the function $\phi^\hbar(z)$ is meromorphic with poles at the
points $\{\pi i(m+ n\hbar)|m,n \in {\bf N}\}$
and $\{-\pi i (m+n\hbar)|m,n \in {\bf N}\}$.
\end{prop}

{\bf Proofs.}
To prove property {\bf P.1} note that $\lim_{z \rightarrow
-\infty}\phi^\hbar(z)=0$.
Therefore, it suffices to prove, that
$$
\lim_{\hbar \rightarrow
0}\frac{\partial}{\partial z} \phi^\hbar(z) = \frac{1}{e^{-z}+1}.
$$
The left-hand side of this equality can be
easily computed using residues:
$$
\lim_{\hbar \rightarrow 0}\frac{\partial}{\partial z} \phi^\hbar(z) =
\lim_{\hbar \rightarrow 0}-\frac{\pi\hbar}{2}\int_\Omega
\frac{-ipe^{-ipz}}{\sinh(\pi
p)\sinh(\pi \hbar p)}dp =
$$
$$
=\frac{i}{2} \int_\Omega\frac{e^{-ipz}}{\sinh(\pi p)}=
\frac{i}{2} \frac{1}{1+e^{z}}\left(\int_{\Omega}
-\int_{\Omega+i}\right)\frac{e^{-ipz}}{\sinh(\pi p)}dp =
\frac{-1}{1+e^{z}}{\mbox{Res}}_{p=i \pi} \frac{e^{-ipz}}{\sinh(\pi p)} =
$$
$$
=\frac{1}{e^{-z}+1}. \eqno{\bf q.e.d.}
$$

Property {\bf P.2} can be verified by computing the right-hand side using
residues:
$$
\phi^\hbar(z)-\phi^\hbar(-z) =
$$
$$
=-\frac{\pi \hbar}{2}\int_\Omega\frac{e^{-ipz}-e^{ipz}}{\sinh(\pi
p)\sinh(\pi \hbar p)}dp=
-\frac{\pi \hbar}{2}\left(\int_\Omega+\int_{-\Omega}\right)
\frac{e^{-ipz}}{\sinh(\pi
p)\sinh(\pi \hbar
p)}dp =
$$
$$
= \frac{\pi \hbar}{2} 2\pi i\,
\mbox{Res}_{z=0}\frac{e^{-ipz}}{\sinh(\pi p)\sinh(\pi \hbar p)}= z.
\eqno{\bf q.e.d.}
$$

Property {\bf P.3} can be obtained by changing the integration
variable $q=-\overline{p}$:
$$
\overline{\phi^\hbar(z)}=-\frac{\pi\hbar}{2}\int_{\overline{\Omega}}
\frac{e^{ip\overline{z}}}{\sinh(\pi p)\sinh(\pi \hbar p)}dp =
$$
$$
= \frac{\pi\hbar}{2}\int_{-\Omega} \frac{e^{ip\overline{z}}}{\sinh(\pi
p)\sinh(\pi \hbar
p)}dp=-\frac{\pi\hbar}{2}\int_{\Omega} \frac{e^{-ip\overline{z}}}{\sinh(\pi
p)\sinh(\pi \hbar
p)}dp=\phi^\hbar(\overline{z}).
\eqno{\bf q.e.d.}
$$

Property {\bf P.4} can be obtained by changing the integration
variable $q = p/\hbar$:
$$
\phi^{1/\hbar}(z/\hbar) =
-\frac{\pi}{2\hbar}\int_{\Omega}\frac{e^{-ipz/\hbar}}{\sinh(\pi p)\sinh(\pi
p/\hbar)}dp=
$$
$$
=-\frac{\pi}{2\hbar}\int_{\Omega}\frac{e^{-iqz}}{\sinh(\pi \hbar q)\sinh(\pi
q)}d(q \hbar)=
\phi^{\hbar}(z)/\hbar. \eqno{\bf q.e.d.}
$$

Properties {\bf P.5} and {\bf P.6} can be proven
analogously to the property {\bf P.1}.
Both proofs are analogous, and we present only the first of them here
$$
\phi^\hbar(z+i\pi\hbar)-\phi^\hbar(z-i\pi\hbar)
= -\frac{\pi\hbar}{2}\int_\Omega \frac{e^{-ipz}(e^{\pi\hbar p} -
e^{-\pi\hbar p})}{\sinh(\pi
p)\sinh(\pi \hbar p)}dp =
$$
$$
=-\pi\hbar\int_\Omega \frac{e^{-ipz}}{\sinh(\pi p)}dp = \frac{-\pi
\hbar}{e^z+1}
\mbox{~Res}_{p=i\pi} \frac{e^{-ipz}}{\sinh(\pi p)} = \frac{2 \pi i
\hbar}{e^{-z}+1}. \eqno{\bf q.e.d.}
$$

Property {\bf P.7} follows from the obvious observation
that integral (\ref{phi}) converges at
$|\Im{z}| < \pi(1+\hbar)$ and, using properties {\bf P.5} and {\bf P.6},
the function
$\phi^\hbar$ can be continued to the whole complex plane.

\subsection{Proof of Theorem \ref{th4} and its corollaries.}

Corollary {\bf C.1} follows from identity {\bf P.1}.

For proving property {\bf C.2} we must verify that the
morphism ${\cal
T}^\hbar(\Gamma)\rightarrow {\cal T}^{1/\hbar}(\Gamma)$ commutes with a flip
(this morphism obviously
commutes with a graph symmetry). This implies that $(A+\phi^\hbar(Z))/\hbar =
A/\hbar +\phi^\hbar(Z/\hbar)$, \
$(B-\phi^{1/\hbar}(-Z))/\hbar = A/\hbar
-\phi^\hbar(-Z/\hbar)$, etc.
Therefore, it suffices to prove that $\phi^\hbar(z)/\hbar =
\phi^{1/\hbar}(z/\hbar)$,
which is just property {\bf P.4}.

For proving the corollary {\bf C.3} we must verify that a flip maps
central generators to the corresponding central generators. This follows
from property {\bf P.2} of the function $\phi^\hbar(z)$.

{\bf Proof of Lemma~\ref{lem1}}.

One can reformulate Lemma~\ref{lem1} as follows. Construct a
sequence of elements of the quantum torus using
the recursion relation $U_{-1} =
V^{-1}$, $U_0 = U$, $U_{i+1}= (1+qU_{i})U_{i-1}^{-1}$. Then, $qU_{i+1}U_i =
q^{-1}U_iU_{i+1}$ and $U_{i+5} = U_i$ for any $i$.

The lemma can be proved by a straightforward calculation:
$$
U_{i+5} = (1+qU_{i+4})U^{-1}_{i+3} = U^{-1}_{i+3} + q U^{-1}_{i+2}
U^{-1}_{i+3} +
U^{-1}_{i+2} =
$$
$$
=(1+q^{-1}U_{i+1})U^{-1}_{i+2} = U_{i}. \eqno{\bf q.e.d.}
$$

The proofs of properties {\bf C.3}, {\bf C.4}, and {\bf C.5}
are straightforward calculations;
remark only that {\bf C.3} and  {\bf C.4} follow directly from
the Stone--von~Neumann theorem.

The proof of Theorem~\ref{th4} can be reduced to proving the following
set of statements.

\vspace{2mm}
{\bf S.1} Morphisms (\ref{q-mor}) are indeed the algebra morphisms.

{\bf S.2} Morphism (\ref{q-mor}) preserve the $*$-structure.

{\bf S.3} Morphisms (\ref{q-mor}) satisfy relations {\bf R.2}, {\bf R.4},
and the graph symmetry relations.

{\bf S.4} Morphisms (\ref{q-mor}) satisfy the pentagon identity {\bf R.5}.
\vspace{2mm}

Statement {\bf S.1} follows from the equalities $[A+\phi^\hbar(Z),
B-\phi^\hbar(-Z)] = 0$, \ $[A+\phi^\hbar(Z), D -
\phi^\hbar(-Z)]= 2\pi i \hbar$ (the remaining relations are analogous or
obvious).
Using the property {\bf P.2} we can transform the commutators:
$$
[A+\phi^\hbar(Z), B-\phi^\hbar(-Z)] =
[A,B]-[A,\phi^\hbar(-Z)]+[\phi^\hbar(Z),B]=
$$
$$
=2\pi i \hbar\left(-1 - \frac{\partial}{\partial Z}(-\phi^\hbar(-Z)) +
\frac{\partial}{\partial Z}
\phi^\hbar(Z)\right)=0.
$$
Analogously,
$$
[A+\phi^\hbar(Z), D - \phi^\hbar(-Z)]=[A,D] -
[A,\phi^\hbar(-Z)]+[\phi^\hbar(Z),D]=
$$
$$
=2\pi i\hbar\left(\frac{\partial}{\partial Z}(-\phi^\hbar(-Z))+
\frac{\partial}{\partial Z}\phi^\hbar(Z)\right)= 2\pi i\hbar.
$$

The proof of the statement
{\bf S.2}, i.e., that morphisms (\ref{q-mor}) preserve the
real structure is obviously equivalent to the
realness condition for
the function $\phi^\hbar(z)$, i.e., to the property {\bf P.3}.

Morphisms (\ref{q-mor}) agree with relation {\bf R.2} because
of property {\bf P.2}. The agreement with relation {\bf R.4} as well as
with the graph symmetry relations is obvious by construction.  The only
nontrivial verification is the proof that morphisms (\ref{q-mor})
reproduce pentagon identity {\bf R.5}

There are seven generators involved in the sequence of flips ${\bf R.5}$.
Denote them by $A_0$, $B_0$, $C_0$, $D_0$, $E_0$,
$X_0$, and $Y_0$ as shown in Fig.~1. Note that the flip results in
the cyclical rotation of the piece of graph shown in Fig.~1.
Denote by $A_i$, $B_i$, $C_i$, $D_i$, $E_i$, $X_i$, and
$Y_i$ the algebra elements associated to the edges of this piece of graph
after performing $i$~flips.
The rules how these elements are changed by (\ref{q-mor}) are
\be
\begin{array}{rcl}
X_{i+1} &=& Y_i - \phi^\hbar(-X_i)\\
Y_{i+1}&=&-X_i\\
A_{i+1}&=&D_i\\
B_{i+1}&=&E_i\\
C_{i+1}&=&A_i + \phi^\hbar(X_i)\\
D_{i+1}&=&B_i-\phi^\hbar(-X_i)\\
E_{i+1}&=&C_i + \phi^\hbar(X_i)
\end{array}
\ee

Our aim is to prove that these sequences of operators are five-periodic.

Assume for a moment that the five-periodicity of $X_i$ is proved. Then the
five-periodicity of $Y_i$
is obvious, because $Y_{i+1}=-X_i$.
The five-periodicity of, say, $A_i$ follows
from the calculation
$$
X_{i+1}=Y_i - \phi^\hbar(-X_i)= -X_{i-1} -\phi^\hbar(-X_i).
$$
Therefore,
$$
\phi^\hbar(-X_i) = -X_{i+1}-X_{i-1}.
$$
Taking into account {\bf P.2}, we have
$$
\phi^\hbar(X_i) = X_i-X_{i+1}-X_{i-1}
$$
Now we can use these identities to transform $A_{i+5}$:
\bea
A_{i+5}&=&D_{i+4}=B_{i+3}-\phi^\hbar(-X_{i+3})=E_{i+2}-\phi^\hbar(-X_{i+3})=
\nonumber\\
&=&C_{i+1}+\phi^\hbar(X_{i+1})-\phi^\hbar(-X_{i+3})=A_i
+\phi^\hbar(X_{i})+\phi^\hbar(X_{i+1})-\phi^\hbar(-X_{i+3})=
\nonumber\\
&=& A_i + (X_i - X_{i-1}-X_{i+1}) + (X_{i+1} - X_{i}-X_{i+2})+(X_{i+4} +
X_{i+2})
\nonumber\\
&=& A_i +
X_{i+4}-X_{i-1}=A_i.
\nonumber
\eea
We have shown that the five-periodicity of $A_i$ (and therefore of $B_i,
C_i,D_i,E_i$, and $Y_i$)
follows from the five-periodicity of $X_i$.

For proving the five-periodicity of $X_i$, we introduce the
algebraic elements
$$
U_i = e^{-X_i}, \qquad \widetilde{U}_i=e^{-X_i/\hbar},
$$
which satisfy the following commutation relations:
\be
qU_iU_{i+1}=q^{-1}U_{i+1}U_i;
~\tilde{q}\widetilde{U}_i\widetilde{U}_{i+1}
=\tilde{q}^{-1}\widetilde{U}_{i+1}\widetilde{U}_i;~
U_i\widetilde{U}_j=\widetilde{U}_jU_i;~\widetilde{U_i} = 
(U_i)^\hbar,\label{comm}
\ee
where
$$
q=e^{\pi i \hbar},\qquad \tilde{q} =e^{\pi i/\hbar}.
$$
These algebraic elements transform in an especially simple way:
\begin{eqnarray}
U_{i+1}&=& (1+qU_{i})U_{i-1}^{-1}\label{e}\\
\widetilde{U}_{i+1}&=&(1+\tilde{q}
\widetilde{U}_{i})\widetilde{U}_{i-1}^{-1}\label{te}
\end{eqnarray}
Indeed,
$$
U_{i+1}=e^{-X_{i+1}}=e^{-Y_i+\phi^\hbar(-X_i)}=\exp\left(\frac{1}{2\pi i
\hbar}\int_{-X_i}^{-X_i+2\pi i \hbar}\phi^\hbar(z)dz\right)e^{-Y_i}=
$$
$$
=\exp\left(\frac{1}{2\pi i \hbar}\int_{-\infty}^{-X_i}(\phi^\hbar(z+2\pi i
\hbar)
-\phi^\hbar(z))dz\right) e^{X_{i-1}}=
$$
$$
= \exp\left(\int_{-\infty}^{X_i}\frac{dz}{e^{-z-\pi i \hbar} +1}\right)
U_{i-1}^{-1} =
(1+qU_i)U_{i-1}^{-1},
$$
where we have used the standard formula
$$
e^{A+\Phi(B)} =
\exp\left\{\frac{1}{[A,B]}\int_B^{B+[A,B]}\Phi(z)dz\right\}e^A,
$$
which is valid for all~$A$ and~$B$ such that the commutator $[A,B]$ is a
nonzero scalar.

The proof of (\ref{te}) is analogous.

Now in order to prove that $X_i$ is five-periodic it suffices to verify,
that both $U_i$ and $\widetilde{U}_i$ are five-periodic. Indeed,
the five-periodicity of $U_i$ only does not suffice because the logarithm
of an operator is ambiguously
defined. However, if we have two families of operators $U$
and $\widetilde{U}$
depending on $\hbar$ continuously, then if there exists an operator $X$
(also depending on $\hbar$ continuously) such that $U=e^X$ and
$\widetilde{U}=e^{X/\hbar}$, then this operator is
uniquely determined. It can be found as the limit
$$
\lim_{(m+n/\hbar)\rightarrow 0} (U^m\widetilde{U}^n)/(m+n/\hbar)
$$
for any irrational value of $\hbar$.

Relations (\ref{e}) and (\ref{te}) coincide with the relations from the
assertion of Lemma~\ref{lem1}, and the sequences $U_i$ and
$\widetilde{U}_i$
are therefore five-periodic.\hfill{\bf q.e.d.}

To prove relation {\bf C.7}, we consider the representation of the
algebra generated by two real generators $X$ and $Y$
with the relation $[X,Y] =
-2\pi i \hbar$. Any unitary representation of this algebra is equivalent to
the representation on the space $L^2(\RR)$, where $X = 2\pi i \hbar
{\partial}/{\partial z}$ and $Y = z$, while $z$ is the parameter on the
real line. Now one can easily verify that $K^{-1}XK = Y-\phi^\hbar(-X)$
and  $K^{-1}YK = -X$. As follows from Theorem~\ref{th4}, $K^{-5} X K^5 = X$
and $K^{-5} Y K^5 = Y$, and the operator $K^5$ is
therefore a scalar.\hfill{\bf q.e.d.}

Relation {\bf C.8} can be proved analogously. The quantum torus
$qUV=q^{-1}VU$ (where $q=e^{\pi i m/n}$ and both $m$ and $n$ are odd) has a
center generated by $U^n$ and $V^n$. For each nonzero constants $u$ and $v$
there exists a unique (up to the conjugation) finite-dimensional irreducible
representation for which $U^n=u$ and $V^n=v$. In this case,
$((1+qU)V)^n = (1+u)v$ and the action of the automorphism on the
set of representations has therefore the order 5 (as follows from the
lemma). We choose the explicit matrix realization: \
$U\{u\}^i_j = u^{1/n} \delta^{i+1}_{j}$, \
$V\{v\}^i_j = v^{1/n}q^{2i}\delta^{i}_{j}$, $i,j \in \ZZ/n\ZZ$.
We now can easily verify that $L(u)U\{u\}L^{-1}(u)
\sim (1+qU\{(1+u)v\})V\{u^{-1}\})$ and $L(u)V\{v\}L^{-1}(u) \sim
U\{(1+u)v\}^{-1}$, where the sign
$\sim$ means the equality up to a scalar factor. If we
denote by $L$ the left-hand side of equality (\ref{dld}),
then $L U\{u\} L^{-1} \sim
U\{u\}$ and $L V\{v\} L^{-1} \sim V\{v\}$.
Therefore, $L = \mbox{const}$ because
the representation with the given value of the center is unique. \hfill{\bf
q.e.d.}

\subsection*{Acknowledgments}
One of the authors (V.F.) is very indebted to the Centre de Physique
Theorique in Marseille, where an important part of this paper was written and
discussed and to V.~Ovsienko and O.~Ogievetsky for
the very valuable discussion.

The work was partially financially
supported by the Russian Foundation for Basic
Research (Grant No.~98-01-00327) and by ork of
the program for the Support of Leading Scientific Schools (Grant
No.~96-15-96455).

\end{document}